\documentclass[12pt]{article}
\usepackage{amssymb}
\usepackage{amsthm}
\usepackage{graphicx}
\usepackage[normalem]{ulem}
\newtheorem{teo}{\bf Theorem}[section]
\newtheorem{pro}[teo]{\bf Proposition}
\newtheorem{lem}[teo]{\bf Lemma}
\newcommand{\NN}{\mathbb{N}}
\newcommand{\RR}{\mathbb{R}}
\newcommand{\QED}{\hfill \rule{6pt}{8pt}}
\title{On the Transitivity of Invariant Manifolds of Conservative Flows}
\author{F\'abio Castro and Fernando Oliveira}
\begin{document}
\maketitle
\begin{abstract}
The main result of this work is the following: for volume preserving flows on compact manifolds with the $C^r$ topology, $1 \leqq r \leqq \infty$ , the closure of every invariant manifold of periodic orbits and singularities is a chain transitive set.

We also develop to new local constructions, which surprise by the simplicity of the arguments. One, a local perturbation to change an orbit to a nearby without altering its past. The other is a flow box theorem in the context of volume preserving flows, a result that is well known for Hamiltonians or general flows.
\end{abstract}

\section{Introduction}
In differentiable dynamics, after periodic orbits one of the most important objects are the invariant manifolds.\

Homocliic orbits also play an essential role.
Poincar\'{e} was probably the first one to draw the
trellis of invariant manifolds that arises in the presence of a tranverse
homoclinic orbit, and became astonished with the complexity of the phenomenon.
He conjectured that periodic orbits would be dense for most conservative systems. In his great Prize Memoir in the Acta Mathematica he dealt with the $n$ body problem. He employed the method of analytic continuation and after some work found it to be insufficient. Then he turned his attention to the restricted three body problem and wrote his seminal work M\'ethods nouvelles de la M\'echanique c\'eleste. By the end of his life he gave out his last geometric theorem without a complete proof. It is very interesting to see how Birkhoff describes the history and events in the work of Poincar\'{e} in his article "On the periodic motions of dynamical systems"  in Acta Mathematica 1927.

In the $C^r$ topology, $r>1$, for symplectic or volume  preserving diffeomorphisms it is not known whether periodic orbits are dense, even in the case of the two sphere. This problem is equivalent to the density of homoclinic points in invariant manifolds of hyperbolic periodic points. The general feeling is that they are true. If this is the case, then the closure of each invariant manifold has a dense orbit. This is not known as well, even for surfaces. If we try something a little weaker like chain transitivity, then its possible to prove that for symplectic and volume preserving diffeomorphisms of compact manifolds, generically the closure of invariant manifolds of hyperbolic periodic  points are chain transitive sets \cite{oliveira}.\

These problems all make sense for conservative flows (volume preserving, hamiltonian, geodesic flows, etc), in manifolds of dimension at least $3$.

The main result of this paper is the following:\\

\noindent {\it Let $M$ be a compact manifold  without boundary and dimension greater or equal to 3. Let $\omega$ be a volume form on $M$ and consider the set ${\cal X}^r_\omega(M)$ of volume preserving vector fields on $M$, with the $C^r$ topology. Then,  ${\cal X}^r_\omega(M)$ contains a residual subset $\cal R$ such that any $X\in {\cal R}$ has the following property: the closure of an invariant manifold of a hyperbolic critical element of X is a chain transitive set $($by a critical element of a vector field, we mean a periodic orbit or a singularity$)$.}\\

\begin{figure}[htbp]
	\centering
\includegraphics[width=0.7\textwidth]{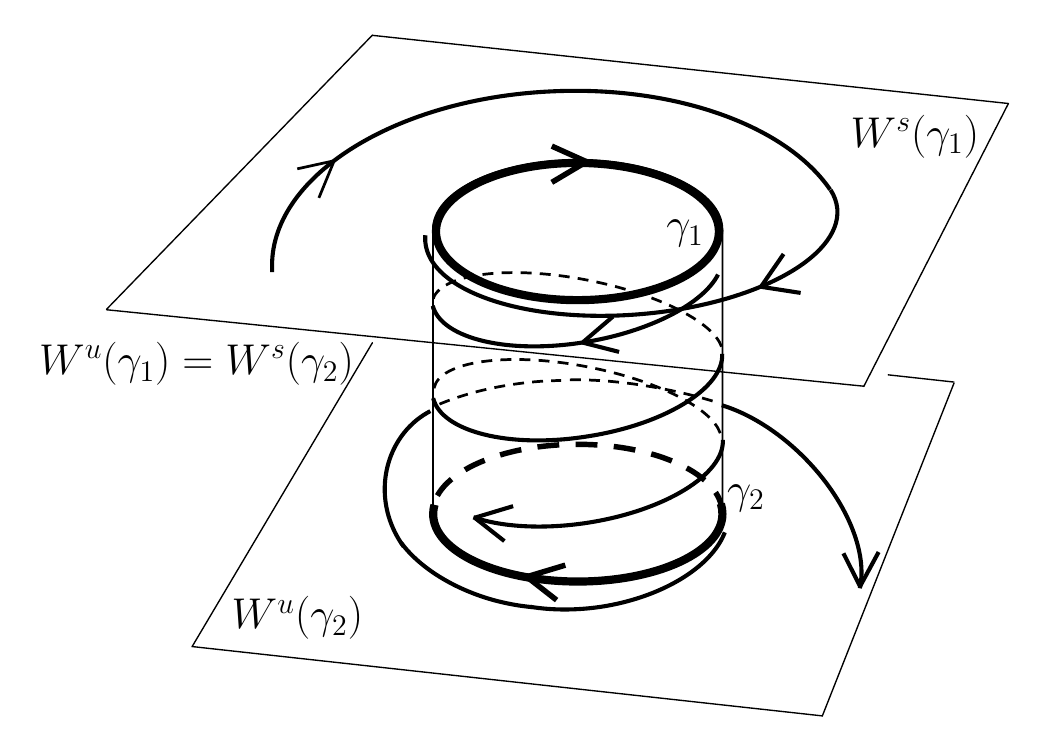}
		\caption{$\overline{W^u(\gamma_1)}=\overline{W^s(\gamma_2)}$ has $\gamma_1$ as a repeller and $\gamma_2$ as an attractor.}
		\label{F2}
\end{figure}

A conservative system obviously does not have attractors or repellers. But being lower dimensional objects in general, closures of invariant manifolds may have attractors and repellers (see figure \ref{F2}). In fact some times they do, but this is not generic. Chain transitive sets have no attractors or repellers.

The proof follows from the fact that generically invariant manifolds contain a dense subset of recurrent orbits.

\section{Some basic facts}
Now we recall some basic facts. Good references are the books of Abraham and Marsden \cite{abraham}, Conley \cite{conley}, Palis and de Melo \cite{palis} and Ma\~n\'e \cite{mane}.\par

We consider a $C^\infty$ compact manifold $M$  without boundary of dimension $n\geq  3$, with a volume form $\omega$ which we also assume to be $C^\infty$. We also consider a Riemannian metric compatible with its topology.\par

Let ${\cal X}^r_\omega(M)$ be the set of $C^r$ vector fields $X$ on $M$ which preserves the volume form $\omega$, $1\leq r\leq \infty$. Equipped with the $C^r$ topology, ${\cal X}^r_\omega(M)$ is a Baire space, and we say that a property is $C^r$ generic if it holds for vector fields in a residual subset of ${\cal X}^r_\omega(M)$. \par 

An $(\epsilon,t)$-chain from $p$ to $q$ is a finite sequence $[x_0,x_1, \dots, x_m]$ of points in $M$ such that $x_0=p$, $x_n=q$ and for each $i\in\{0,\dots,m-1\}$ there exists $t_i\geq t$ such that $d(X_{t_i}(x_i), x_{i+1})<\epsilon$.\par
An invariant subset $A$ by the flow $X_t$ is {\it chain transitive} if there exists an $(\epsilon,t)$-chain in $A$ from $p$ to $q$, for every  $p,q\in A$, $\epsilon>0$ and $t>0$.

Let $\gamma$ be a periodic orbit of $X$ and $p\in \gamma$. Let $\Sigma$ be a section transversal to $X$ at $p$ and $\Sigma_p$ a neighborhood of $p$ in $\Sigma$, where the Poincaré map $P:\Sigma_p\to \Sigma$ is defined. We say that $\gamma$ is hyperbolic if $p$ is a hyperbolic fixed point of $P$. 

The unstable manifold of $\gamma$ is the set $W^u(\gamma)=\{x\in M|\; \alpha(x)=\gamma\}$, where $\alpha(x)$ denotes the $\alpha$-limit of the orbit of $x$. Similarly, the  stable manifold of $\gamma$ is the set $W^s(\gamma)=\{x\in M|\; \omega(x)=\gamma\}$, where $\omega(x)$  denotes the $\omega$-limit of the orbit of $x$. The sets $W^u(\gamma)$ and $W^s(\gamma)$ are immersed connected submanifolds of $M$. For $p\in \gamma$, let $\Sigma_p^u$ and $\Sigma_p^s$ be the intersection of $\Sigma_p$ with the local invariant manifolds of $\gamma$.\par

The sets $D^u(p)=P(\Sigma^u_p)-\Sigma^u_p$ and $D^s(p)=\Sigma^s_p-P(\Sigma^s_p)$ are called unstable and stable fundamental domains of $W^u(\gamma)$ and $W^s(\gamma)$, respectively. We have $W^u(\gamma)-\gamma=\bigcup_{t\in \RR}X_t(D^u(p))$ and  $W^s(\gamma)-\gamma=\bigcup_{t\in \RR}X_t(D^s(p))$.

Let $a\in M$ be a hyperbolic singularity of $X$, that is, $DX(a):TM_a\to TM_a$ does not have eigenvalues on the imaginary axis. The unstable manifold of $a$ is the set $W^u(a)=\{x\in M|\; \alpha(x)=a\}$. Similarly, the  stable manifold of $a$ is the set $W^s(a)=\{x\in M|\; \omega(x)=a\}$. Let $S$ be a small sphere centered at $a$. The sets $S^u(a)$ and $S^s(a)$  obtained by the intersection of $S$ with the local unstable and stable  invariant manifolds of $a$, are called fundamental domains of $W^u(a)$ and $W^s(a)$. We have that $W^u(a)-a=\bigcup_{t\in \RR} X_t(S^u(a))$ and $W^s(a)-a=\bigcup_{t\in \RR} X_t(S^s(a))$.

\section{Tubular flow for volume preserving vector fields}

Now, we describe a result (flow box) which is very well known for general vector fields and Hamiltonians, but not for volume preserving flows. See Herman \cite{herman} for a nice complete treatment, and Bessa \cite{bessa} for a three dimensional version based on different techniques.
\begin{lem}[Tubular flow for volume preserving vector fields] Let $X\in \mathcal{X}^\infty_\omega(M)$, $a\in M$ be a regular point of the vector field $X$ and let $\Sigma$ a transversal section to $X$ which contains $a$.  Then there exists a $C^{\infty}$ coordinate system $\alpha:U\subset M\to V\subset \RR^n$ with $\alpha(a)=0$ and such that: 
\begin{enumerate}
	\item $\alpha_*X\equiv(0,\dots,0,1)$.
	\item $\alpha_*\omega= dx_1\wedge\cdots\wedge dx_n$.
	\item $\alpha(U\cap \Sigma)\subset \{x_n=0\}$.
\end{enumerate}
 \label{lema:vol}
\end{lem}
{\it Proof.} The known Tubular Flow Theorem gives a $C^\infty$ chart $\varphi:V_1\to M$  with $\varphi(0)=a$ and such that $\varphi^* X(z_1,\dots,z_n)=(0,\dots,0,1)$, for all $z=(z_1,\dots,z_n)\in V_1$. Taking the {\it pull-back} form $\varphi^*\omega$, we have $\varphi^*\omega(z)=\psi(z)dz_1\wedge \cdots \wedge dz_n$, where $\psi:V_1\to \RR$ is a $C^\infty$ function, which we may assume to be positive by taking a change of orientation in $M$, if necessary. \par 

Denote $\varphi^*X=Z=(Z_1,\dots,Z_n)\equiv (0,...,0,1)$ and $\varphi^*\omega=\eta$. Then it follows that 
  $$\displaystyle div_\eta Z(z)=\sum_{i=1}^n \frac{\partial Z_i}{\partial z_i}(z)+ \frac{1}{\psi(z) }\sum_{i=1}^n \frac{\partial \psi}{\partial z_i}(z)Z_i(z),$$ 
	
	where $div_\eta Z$ is the function which satisfies Liouville's formula, that is, $$\mbox{det}_\eta\phi'_t(z)=\exp \int_0^t(\mbox{div}_\eta Z)(\phi_s(z))ds,$$ where $\phi_t$ is the time $t$ map of the flow of $Z$ (for details see page 130 of Abraham and Marsden \cite{abraham}). 
	
	Since $X$ preserves $\omega$, we have that $Z$ preserves $\eta$. Note that 
	$$ \sum_{i=1}^n \frac{\partial Z_i}{\partial z_i}(z)+ \frac{1}{\psi(z)}\sum_{i=1}^n \frac{\partial \psi}{\partial z_i}(z)Z_i(z)=\frac{\partial \psi}{\partial z_n}(z).$$ 
Therefore, by the Liouville's formula, we have
$$\frac{\partial \psi}{\partial z_n}(z)=0,$$  
for all $z\in V_1$, that is, $\psi$ does not depend on $z_n$. Now, let $f(z_1,\dots,z_{n-1})=\psi(z_1,\dots,z_{n-1},z_n)$, and let $\xi:V_1\to \RR^n$ defined by $$\xi(z_1,\dots,z_n)=(z_1,\dots,z_{n-2},\int_0^{z_{n-1}}f(z_1,\dots,z_{n-2},t)dt, z_n),$$ 
whose jacobian is of the form
$$J\xi(z)=\left[
\begin{tabular}{ccccccc}
1&0&$\cdots$&0&0&0&0\\
0&1&$\cdots$&0&0&0&0\\
$\vdots$&$\vdots$&$\ddots$&\vdots&$\vdots$&$\vdots$&\vdots\\
0&0&$\cdots$&1&0&0&0\\
0&0&$\cdots$&0&1&0&0\\
$\ast$&$\ast$&$\cdots$&$\ast$&$\ast$&$f(z')$&0\\
0&0&$\cdots$&0&0&0&1
\end{tabular}\right],
$$

\noindent where $z'=(z_1,\dots,z_{n-1})$. Then, $\xi_*Z=Z$. \par

Note that $\xi(0)=0$ and $\det(J\xi(0))=f(0)=\psi(0,\dots,0,z_n)>0$, since $\psi$ is positive on $V_1$. Therefore, there exist neighborhoods of $0$, $V_2\subset V_1$ and $V_3$, such that the restriction $\xi:V_2\to V_3$ is a diffeomorphism.

\par
On the other hand, if $y_1,\dots,y_n$ are the coordinate functions of $\xi$, that is $\xi(z_1,\dots,z_n)=(y_1,\dots,y_n)$, we have
$$
\left\{\begin{tabular}{l}
$dy_i=dz_i$, if $i\not=(n-1)$\\
$dy_{n-1}=fdz_{n-1}+[\mbox{terms with}\;dz_j], \;j\not=(n-1).$
\end{tabular}\right.
$$

From this, \\
$dy_1\wedge \cdots \wedge dy_n=$\\
$dz_1\wedge \cdots \wedge dz_{n-2}\wedge\Big\{fdz_{n-1}+[\mbox{terms with}\;dz_j],\;j\not=(n-1)\Big\}\wedge dz_n=\\
fdz_1\wedge \cdots\wedge dz_n=\psi dz_1\wedge \cdots\wedge dz_n=\eta$.\\
 This proves that $\xi^*(dy_1\wedge \cdots \wedge dy_n)=\psi dz_1\wedge \cdots\wedge dz_n=\phi^*\omega$.

Now, taking $U=\varphi(V_2)$ and $\beta:U\to V_3$ given by $\beta=\xi\circ\varphi^{-1} $,  $\beta$ satisfies itens 1 and 2.\par

For 3, consider $U$ small so that $\Sigma_0=U\cap \Sigma$ is a connected set that contains  $a$. In this way,  $\beta(\Sigma_0)$ is the graph of a function  $g:A\subset\RR^{n-1}\to \RR$. Now take $\gamma:V_3\to \RR^n$ given by
$$\gamma(y_1,\dots,y_n)=\Big(y_1,\dots,y_{n-1},y_n-g(y_1,\dots,y_{n-1})\Big).$$
\par $\gamma $ is a translation along each line of the flow. It is easy to see that  $\gamma$ preserves the volume form, the vector field and takes the graph of $g$ onto a subset of $\{x_n=0\}$. Finally, we take $V=\gamma(V_3)$ and $\alpha:U\to V$ given by $\alpha=\gamma\circ\xi\circ\phi^{-1}$. $\alpha$ is the desired function.  \QED

\section{Perturbation Lemmas}

Let $B_\delta(x,y)\subset \RR^{2}$ be the open ball of center $(x,y)$ and radius $\delta$. Similarly, $B_\delta[x,y]$ will denote the closed ball. If $C$ is the cylinder $\partial B_\delta(x,y)\times [0,h]\subset \RR^{3}$ and $0<\xi<\delta$,  we define a neighborhood of $C$ as 
$$
A_\xi(C)=\Big(B_{\delta+\xi}(x,y)-B_{\delta-\xi}[x,y]\Big)\times [0,h]\subset \RR^{3}. 			
$$
and call it {\it cylinder ring} with center at $C$ and radius $\xi$.
\begin{lem}[Local perturbation in $\RR^3$]
Let $X:\RR^3\to \RR^{3}$ be the constant vector field defined  by $X(x,y,z)=(0,0,1)$. Consider the cylinder $C= \partial B_\delta(0,0)\times [0,h]\subset \RR^{3}$, $\delta>0$, $h>0$, and points $p\in \partial B_\delta(0,0)\times \{0\}$ and $q\in \partial B_\delta(0,0)\times \{h\}$. Let $\theta$ be the angle between the vectors $p-(0,0,0)$ and $q-(0,0,h)$. Given $\xi>0$  $(\xi<\delta)$,  there exist a $C^\infty$ vector field $Z$ on $\RR^{3}$ with the following properties:
  
\begin{enumerate}
	\item $Z$ preserves the canonical volume form $dx\wedge dy\wedge dz$;
	\item $Z\equiv X$ outside the cylinder ring $A_\xi(C)$;
	\item The positive orbit of $p$, with respect to $Z$, contains $q$;
	\item Given $r\in \NN$ and $\epsilon>0$,  if $|\theta|$ is small enough, then $\|Z-X\|_r<\epsilon$, where $\|.\|_r$ denotes the $C^r$ norm on the set of $C^r$ vector fields.
\end{enumerate}
\label{lema:perturb2}
\end{lem}
 
{\it Proof.} By performing a rotation in the coordinates $x$ and $y$, we may assume $p=(\delta,0,0)$ and $q=(\delta\cos \theta, \delta\sin\theta,h)$,	with $-\pi<\theta\leq\pi$. Now we are going to choose a couple of $C^\infty$ bump functions\par

	$\lambda:\RR\to [0,+\infty)$ with support in $[0,h]$ and $\int_{0}^h\lambda(s)ds=1$ and \par
	$\gamma:[0,+\infty)\to [0,+\infty)$ with support in $[\delta-\xi,\delta+\xi]$ and $\gamma(\delta)=1$.\\

\begin{figure}[htbp]
	\centering
\includegraphics[width=0.7\textwidth]{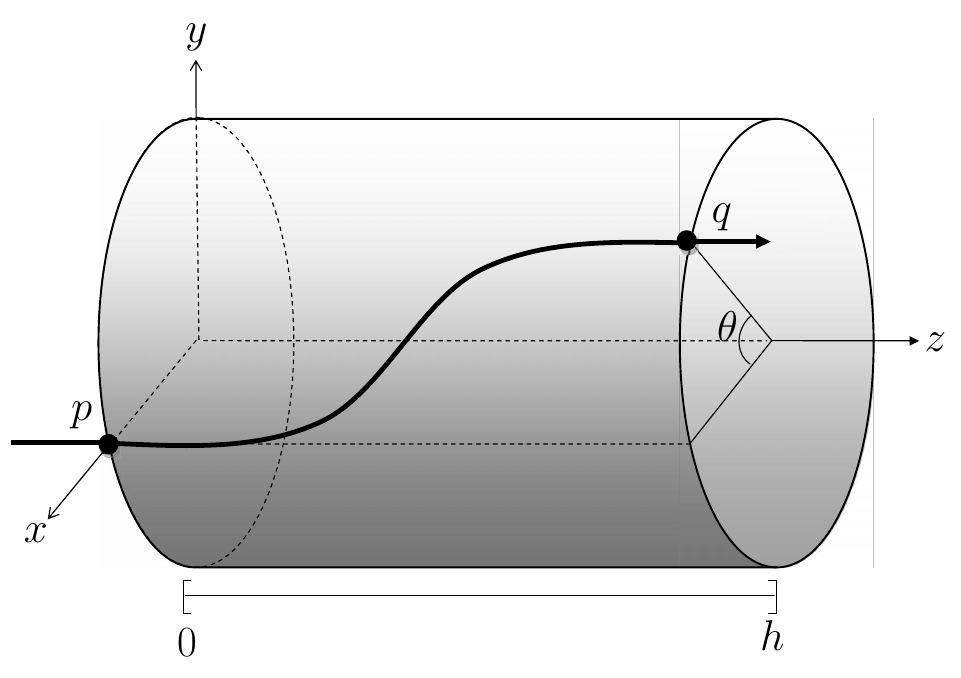}
		\caption{Deviation of the orbit.}
		\label{F1}
\end{figure}

  Consider the vector field $Z:\RR^{3}\to \RR^{3}$ given by

 \begin{equation}
Z(x,y,z)=\left(-yk \lambda(z)\gamma\Big(\sqrt{x^2+y^2}\Big),x k\lambda(z)\gamma\Big(\sqrt{x^2+y^2}\Big) ,1\right),	
	\label{y2}
\end{equation}
where  $k$ is a constant that will be chosen later. \par

Since $\gamma$ is zero on a neighborhood of the origin, the coordinates of $Z$ are $C^{\infty}$ functions, and so is $Z$. It is easy to see that the divergence of $Z$ with respect to the canonical volume form is zero, and $Z$ satisfies item 1. Since $\lambda=0$  outside $(0,h)$ and $\gamma=0$ outside $(\delta-\xi,\delta+\xi)$, $Z=Y$ outside $A_\xi(C)$ and $Z$ satisfies item 2. 

Now we are going to prove item 3. 
We are going ``flow a rotation on the $x$ and $y$ coordinates inside the cylinder ring along the $z$ axis''. This can be seen through the difference of $Z$ and $X$
$$Z-X=\left(-yk \lambda(z)\gamma\Big(\sqrt{x^2+y^2}\Big),x k\lambda(z)\gamma\Big(\sqrt{x^2+y^2}\Big) ,0\right)$$
$$=k\lambda(z)\gamma\Big(\sqrt{x^2+y^2}\Big)\left(-y,x,0\right),$$

``where the classical rotation $(x,y)\mapsto(-y,x)$ appears'' (see figure \ref{F1}).

 The flow of $Z$ is: \par
if $x^2+y^2= 0$, $\varphi_t(0,0,z)=(0,0,t+z)$ and,\par
if $x^2+y^2\not= 0$, then $\varphi_t(x,y,z)=(f(t,x,y,z),g(t,x,y,z),h(t,x,y,z)),$ 
where \\
\begin{tabular}{rl}
	$f(t,x,y,z)=$&$\sqrt{x^2+y^2}\cos\Big(k\gamma(\sqrt{x^2+y^2})\int_z^{t+z}\lambda(s)ds+\arccos \frac{x}{\sqrt{x^2+y^2}}\Big)$, \\
	$g(t,x,y,z)=$&$\sqrt{x^2+y^2}\sin\Big(k\gamma(\sqrt{x^2+y^2})\int_z^{t+z}\lambda(s)ds+\arcsin \frac{y}{\sqrt{x^2+y^2}}\Big)$, \\
	$h(t,x,y,z)=$&$t+z$.  
\end{tabular}\\
Thus, taking $t=h$ and $\displaystyle k=\theta$, we have
\begin{eqnarray}\nonumber
	\varphi_{h}(p)&=&\varphi_{h}(\delta,0,0)\\
	      &=&\Big(\delta\cos\Big[ k\cdot\gamma(\delta)\int_{0}^h\!\!\!\lambda(s)ds\Big],\delta\sin\Big[ k\cdot\gamma(\delta)\int_{0}^h\!\!\!\lambda(s)ds\Big],h\Big) \nonumber \\
	&=&(\delta\cos \theta, \delta\sin \theta, h)=q. \nonumber 
\end{eqnarray}
This concludes the proof of item 3.\par

We can make $\|Z-X\|_r$ as small as we want, by taking $\theta$ as small as necessary. This concludes the proof.\QED\\


Now we would like to see the perturbation result in higher dimension. \par

Firstly, a little notation.  Let $B_\delta(x_1,\dots,x_n,y)\subset \RR^{n+1}$ be the open ball of center $(x_1,\dots,x_n,y)$ and radius $\delta$, and $B_\delta[x_1,\dots,x_n,y]$  the closed ball. If $C$ is the cylinder $\partial B_\delta(x_1,\dots,x_n,y)\times [0,h]\subset \RR^{n+2}$ and $0<\xi<\delta$,  we define a neighborhood of $C$ as 
$A_\xi(C)=\Big(B_{\delta+\xi}(x_1,\dots,x_n,y)-B_{\delta-\xi}[x_1,\dots,x_n,y]\Big)\times [0,h]\subset \RR^{n+2}$
and call it {\it cylinder ring} with center at $C$ and radius $\xi$.\par

\begin{lem}[Local perturbation in higher dimension] Let $X:\RR^n\times\RR\times\RR\to \RR^{n+2}$,  $X(x_1,\dots,x_n,y,z)=(0,\dots,0,0,1)$,  $C= \partial B_\delta(0,\dots,0)\times [0,h]\subset \RR^{n+2}$, $\delta>0$, $h>0$, and points $p\in \partial B_\delta(0,\dots,0)\times \{0\}$ and $q\in \partial B_\delta(0,\dots,0)\times \{h\}$. Let $\theta$ be the angle between the vectors $p-(0,\dots,0,0)$ and $q-(0,\dots,0,h)$. Given $\xi>0$ $(\xi<\delta)$,  there exists a $C^\infty$ vector field $Z$ on $\RR^{n+2}$ with the following properties:\par
1. $Z$ preserves the canonical volume form $dx_1\wedge \cdots \wedge dx_n\wedge dy\wedge dz $;\par
2. $Z\equiv X$ outside the cylinder ring $A_\xi(C)$;\par
3. The positive orbit of $p$, with respect to $Z$, contains $q$;\par
4. Given $r\in \NN$ and $\epsilon>0$,  if $|\theta|$ is small enough, then $\|Z-X\|_r<\epsilon$.
\label{higher}
\end{lem}
{\it Proof.} By performing a rotation in the coordinates $x_n$ and $y$, we may assume
$p=(0,\dots,0,\delta,0,0)$ and $q=(0,\dots,0,\delta\cos \theta, \delta\sin\theta,h)$,	
with $-\pi<\theta\leq\pi$. Consider the vector field $Z:\RR^{n+2}\to \RR^{n+2}$ given by
$$Z(\tilde{x},y,z)=\left(0,\dots,0,-y\theta \lambda(z)\gamma\Big(\sqrt{\|\tilde{x}\|^2+y^2}\Big),x_n \theta\lambda(z)\gamma\Big(\sqrt{\|\tilde{x}\|^2+y^2}\Big) ,1\right),$$
where $\tilde{x}=(x_1,\dots,x_n)$ and $\|\tilde{x}\|^2=x^2_1+\cdots+x_n^2$.
The vector field $Z$ is $C^{\infty}$, the divergence of $Z$ with respect to the canonical volume form is zero and $Z$ satisfies item 1, 2 and 4. Note that the set $Q=\{(0,\dots,0,x_n,y,z)\in \RR^{n+2}|\;x_n,y,z\in \RR\}$
is invariant for $Z$ and $X$, and that $p$ and $q$ belong to $Q$. $Q$ is isomorphic to $\RR^3$ and all calculations will be done in this subspace exactly as in the previous lemma. On the other hand, the support of the perturbation is contained in  $A_\xi(C)$ whose interior is not empty.\par
The flow of $Z$ is \\
$\varphi_t(\tilde{x},0,z)=(\tilde{x},0,t+z)$ if $x_n^2+y^2= 0$;\\ other\-wise,  $\varphi_t(\tilde{x},y,z)=(x_1,\dots,x_{n-1},f(t,\tilde{x},y,z),g(t,\tilde{x},y,z),t+z),$ where \\
$f(t,\tilde{x},y,z)=\sqrt{x_n^2+y^2}\cos\Big(\theta\gamma(\sqrt{\|\tilde{x}\|^2+y^2})\int_z^{t+z}\lambda(s)ds+\arccos \frac{x_n}{\sqrt{x_n^2+y^2}}\Big)$, \\
$g(t,\tilde{x},y,z)=\sqrt{x_n^2+y^2}\sin\Big(\theta\gamma(\sqrt{\|\tilde{x}\|^2+y^2})\int_z^{t+z}\lambda(s)ds+\arcsin \frac{y}{\sqrt{x_n^2+y^2}}\Big)$.\\
 Therefore, $\varphi_{h}(p)=q$ and this concludes the proof of item 3. \QED

\section{Return Lemma}

\begin{lem}  Let $X$ be a $C^\infty$ volume preserving  vector field on $M$. Given $\mathcal{V}$ a $C^r$ neighborhood  of $X$, $p$ a regular point of $M$, $U\subset M$ a neighborhood of  $p$, and $\Sigma$ a transversal section to $X$ that contains $p$, then there exists a $C^\infty$ volume preserving vector field  $K\in \mathcal{V}$ such that:

\begin{enumerate}
	\item (a) $K\equiv X$ outside of   $U$. \par
	 (b) Furthermore, there exists $\tau>0$ such that $K\equiv X$ outside of a subset of $U$ of the form $X_{[0,\tau]}(\Sigma_0):=\{X_t(x)|\;x\in \Sigma_0, 0<t<\tau\}$, where $\Sigma_0$ is a neighborhood of $p$ in $\Sigma$.
	\item Given $T>0$, $K_t(p)\in U$ for some $t>T$.
	\item If $p\notin \alpha(\gamma)$, then $K_t(p)=X_t(p)$, for all $t\leq 0$, where $\gamma$ is the orbit of $p$.
\end{enumerate}

\label{lema-rec1}
\end{lem}

 {\it Proof.} By Lemma \ref{lema:vol}, we may consider a coordinate system $\alpha:U'\to V$ in a neighborhood $U'\subset U$ of $p$, such that $\alpha(p)=(0,\dots,0)$, $\alpha_*X(x_1,\dots,x_{n})=(0,\dots,0,1)$, $\alpha_*\omega=dx_1\wedge\cdots \wedge dx_{n}$, and  $\alpha(\Sigma\cap U')\subset \{x_n=0\}$. We will say that $\alpha(x)$ is recurrent when $x$ is recurrent in $M$.\par
Without loss of generality, we may assume that the coordinates belong to the open set, $V=(-2,2)^{n}$. Let $W=(-1,1)^{n}$, $0<\xi<\delta<1/3$ and $0<h<1$, where $\xi$, $\delta$ and $h$ are the numbers defined in the cylinder ring of Lemmas \ref{lema:perturb2} and \ref{higher}.  Let $\Pi=V\cap \{(x_1,\dots,x_{n})\in \RR^n;\;x_n=0\}$ be a set which is transverse to the flow of $\alpha_*X$. By Poincaré's recurrence theorem, $U'$ contains a dense subset of segments of recurrent points. Therefore, we may take recurrent points in $\Pi$ as close as we want to 0. For every angle $\theta>0$, let us take $q'\in \Pi$ recurrent, whose distance to $0$ is less than $\theta \delta/2$. 
Let $S_{q'}$ be the intersection of $\Pi$ with any subspace of dimension 2  which contains the points $0$ and $q'$. The orbits of $0$ e $q'$ can be seen in $S_{q'}\times (-1,1)$, which is a three dimensional box with two faces parallel to $\{x_n=0\}$. This is the place where we are going to make the perturbation. In $S_{q'}$, consider a circle $C_{q'}$ of radius $\delta$, which contains the points $0$ e $q'$. Observe that the angle $0\hat{O}q'$ is smaller than $\theta$, where $O$ is the center of $C_{q'}$. Let $C_{q'}$ flow to make up the cylinder $C=C_{q'}\times [0,h]$. Let $q$ be the intersection of the orbit of $q'$ with the circle $C_{q'}\times \{h\}$. Now we may apply the perturbation Lemma \ref{higher} at the cylinder $C$ to join $0$ to $q$ and obtain a vector field $Z$ in $W$ such that:\par
	a. $Z$ preserves the canonical volume form;\par
	b. $Z\equiv (0,\dots,0,1)$ in  $S_{q'}\times (-1,1)-A_\xi(C)$;\par
	c. The positive orbit of $0$, with respect to $Z$, contains $q$.\\ 
 Let us define the vector field $K$ in $M$ in the following way: $K\equiv X$ off of $U'$ and $K\equiv\alpha^*Z$ in $U'$.  Note that $K$ is $C^\infty$, satisfies item 1(a), and taking $\theta$ sufficiently small, we may assume $K\in \mathcal{V}$.\par 
In order to proof 1(b), we consider $\Pi_0\subset \Pi$ a compact neighborhood of the origin contained in $\alpha(\Sigma\cap U')\cap W$. Then just take $\Sigma_0=\alpha^{-1}(\Pi_0)$ and $\tau=\sup\{t>0|\;\;\alpha(X_t(x))\in \Pi_0\times[0,h], x\in \Sigma_0\}.$\par
 In order to prove item 3, it is enough to take a neighborhood $U'$ of $p$ sufficiently small, in a way that the negative orbit of $p$ intersects $U'$ in only one connected segment that ends at $p$. This is possible under the hypothesis that $p\notin \alpha(\gamma)$. Hence, the perturbation does not change the negative orbit of $p$. \par
The proof of item 2 is the following.
Consider $\Sigma_0$ and $\tau$ as in item 1(b) and $p_2$ the first return of $\alpha^{-1}(q)$ (which is recurrent) to $\Sigma_0$. If $p_2=p$, there is nothing to prove, otherwise repeat the whole previous argument on the vector field $K$, to join the point $p_2$ to a recurrent point $q_2$ inside the open set $\alpha^{-1}(W)$, without modifying $K$ along its orbit from $p$ to $p_2$. From item 1(b), the new perturbation may be done in the interior of a set  such as $X_{[0,\tau_2]}(\Sigma_1)$ with $p_2\in\Sigma_1\subset \Sigma_0$ and $0<\tau_2<\tau$, with $\tau_2$ and $\Sigma_1$ as small as necessary. This adds at least two units to the return time of $p$ to $\Sigma_0$, since the vector field $\alpha_*K=Z$ is constant of modulus one, and spends time two in $V-W$. Doing this as many times as necessary, we can make a return time as large as wanted. \QED \\

As we have seen, our Perturbation and Return Lemmas are applied only to $C^\infty$ vector fields. The next lemma, due to Carlos Zuppa (\cite{zuppa}), tells us that this is not a problem, because the  $C^\infty$ volume preserving vector fields are dense in ${\cal X}^r_\omega(M)$ under $C^r$ topology. 
 
\begin{lem}[Regularization]
${\cal X}^\infty_\omega(M)$ is dense in ${\cal X}^r_\omega(M)$ in the $C^r$ topology. 
\label{zuppa}
\end{lem}

\section{Recurrence on invariant manifolds} 
\begin{pro}  If dim$(M)\geq 3$ and $1\leq r\leq\infty$, the set $\mathcal{X}^r_\omega(M)$ contains a  residual subset $\mathcal{R}$ such that if $\gamma$ is a hyperbolic critical element of $X\in \mathcal{R}$, then the unstable manifold of $\gamma$ contains a dense subset of $\omega-$recurrent points, and the stable manifold of $\gamma$ contains a dense subset of $\alpha-$recurrent points.
\label{rec}
\end{pro}
\proof Firstly, we would like to observe that it is enough to prove that $\mathcal{X}^r_\omega(M)$ contains a residual subset $R^+$, such that if $\gamma$ is a hyperbolic critical element of $X\in R^+$, then $W^u(\gamma)$ contains a dense subset of $\omega-$recurrent points. This follows from the fact that the map $i(X)=-X$ is a homeomorphism of $\mathcal{X}^r_\omega(M)$, and therefore $R^-=i(R^+)$ is also residual. Let $R=R^+\cap R^-$. If $X\in R$ and $\gamma$ is a hyperbolic critical element of $X$, since $X\in R^+$ then $W^u(\gamma)$ contains a dense set of $\omega-$recurrent points, and since $X\in R^-$, $-X\in R^+$, and $W^s(\gamma,X)=W^u(\gamma,-X)$ contains a dense subset of $\alpha-$recurrent points.\par
Before the Baire category argument, let us make an observation about the construction of the Poincaré map and fundamental domains. The transversal sections $\Sigma$ and $\Sigma_p$ through $p\in \gamma$ may be taken as submanifolds without boundary. So we may assume the same for $\Sigma^u_p$.
The fundamental domain was defined as $D=P(\Sigma^u_p)-\Sigma^u_p$.
We will need to consider sets $D^\circ =P(\Sigma^u_p)-\overline{\Sigma^u_p}$, which are just $D$ minus its boundary as a submanifold of 
 $M$. If $Z$ is a dense subset of $D$ or  $D^\circ$, then $\bigcup_{t\in \RR}(X_t(Z))$ is dense in $W^u(\gamma)$.\par
For each $1\leq r<\infty$, we consider a collection of sets $A^r_{n,k,m}\subset \mathcal{X}^r_\omega(M)$, indexed by $n,k,m\in \NN$, defined as follows: 
$X\in A^r_{n,k,m}$ if it satisfies the following two conditions:

\begin{enumerate}
		\item [(P1).] If $\dim M\geq 4$, then all periodic orbits of $X$ with period  $\leq n$ are hyperbolic. If $\dim M=3$, then all periodic orbits of period  $\leq n$ are elementary; that is, $1$ is not an eigenvalue of the Poincaré map (the orbit is elliptic or hyperbolic).\\

		\item[(P2).] If $\gamma$ is a hyperbolic periodic orbit of $X$ with period $\leq n$, then $W^u(\gamma)$ has a fundamental domain $D$  depending only on $X$,  $\gamma$ and not on $k$ nor $m$, which satisfies the following conditions: 
		
		\begin{itemize}
		\item[(A)] There exist a finite subset $Z(k)$ contained in $D^\circ$ depending only on $X$, $\gamma$ and $k$, such that: 
		\begin{itemize}
			\item[(A1)] $Z(k)$ is $\frac{1}{k}$-dense in $D^\circ$; 
			\item[(A2)] $z\in Z(k)\Rightarrow\;\; X_t(z)\in B_\frac{1}{m}(z),\;\;\mbox{to some $t>m$}.$ 
		\end{itemize}
		\end{itemize}
\end{enumerate}	

The sets $\tilde{A}^r_{k,m}\subset \mathcal{X}^r_\omega(M)$, $k,m\in \NN$, are defined in a similar way: $X\in \tilde{A}^r_{k,m}$ when $X$ satisfies the following conditions:
\begin{enumerate}
		\item [(S1).] All singularities of $X$ are hyperbolic.
		\item [(S2).] If $a$ is a hyperbolic singularity of $X$, then $W^u(a)$ has a fundamental domain $D$ depending only on $X$ and $a$, which satisfies the following conditions:
		\begin{itemize}
				\item[(B)] There exist a finite subset $Z'(k)$ in $D$ depending only on $X$, $a$ and $k$, such that:
		\begin{itemize}
			\item[(B1)] $Z'(k)$ is $\frac{1}{k}$-dense in $D$; 
			\item[(B2)] $z\in Z'(k)\Rightarrow\;\; X_t(z)\in B_\frac{1}{m}(z),\;\;\mbox{for some $t>m$}.$  
		\end{itemize}
		\end{itemize}
\end{enumerate}

We affirm that  $\mathcal{R^+}=\left(\bigcap_{n,k,m}A^r_{n,k,m}\right)\cap \left(\bigcap_{k,m}\tilde{A}^r_{k,m}\right)$ is the set desired. Indeed, we divide the argument in two cases:\par 

1) [periodic orbit] Let $X\in R^+$ and $\gamma$ a hyperbolic periodic orbit of $X$ of period less than or equal to $n$. $X\in A^r_{n,k,m}$ for every $k$ and $m$. Let $D$ be a fundamental domain of $\gamma$ which satisfies (A) for every $k$ and $m$. $Z=\bigcup_{k\geq 1}Z(k)$ is a dense subset of $D$. Let $z\in Z$ and fix $k$ such that $z\in Z(k)$. Since $X\in A^r_{n,k,m}$ for every $m$, $X_t(z)\in B_{1/m}(z)$ for some $t>m$. This implies that $z$ is $\omega-$recurrent. Hence, $\bigcup_{t\in \RR}X_t(Z)$ is a dense subset of $W^u(\gamma)$, made of $\omega-$recurrent orbits.\par

2) [singularity] Let $X\in R^+$ and $a$ a hyperbolic singularity of $X$. We have that $X\in \tilde{A}^r_{k,m}$ for every $k$ and $m$. Let $D$ be a fundamental domain of $a$ which satisfies (B) for every $k$ and $m$. Then $Z'=\bigcup_{k\geq 1}Z'(k)$ is a dense subset of $D$. Let $z\in Z$ and fix $k$ such that $z\in Z'(k)$. Since $X\in \tilde{A}^r_{k,m}$ for every $m$, $X_t(z)\in B_{1/m}(z)$ for some $t>m$, and $z$ is $\omega-$recurrent. Hence $\bigcup_{t\in \RR}X_t(Z')$ is a subset of $W^u(a)$ made of $\omega-$recurrent orbits. \par
 
Note that the proof of the two cases is almost same, except that in the case of a singularity there is no period to deal with. \par

Now let us show that each $A^r_{n,k,m}$ is open.\par
It is well known that the set of vector fields that satisfy condition (P1) forms an open and dense subset of ${\cal X}^r_\omega(M)$. See for example Robinson \cite{robinson1}.  To see that the set of vector fields which satisfy (P2) is open, note that if $X$ satisfies (P1), then the set of hyperbolic periodic orbits of $X$ with period $\leq n$ is finite, has constant cardinality and depends continuously on vector fields in a neighborhood of $X$.\par

Since (A), (A1) and (A2) are open conditions defined in terms of finitely many hyperbolic periodic orbits, finite subsets of fundamental domains of their local invariant manifolds and Poincaré maps, we have that the set of vector fields that satisfy (P2) is open.\par

As before, the sets $\tilde{A}^r_{k,m}$ are open. The arguments are almost the same, except for the period and the topology of fundamental domains. \par

It remains to prove that $A^r_{n,k,m}$ and $\tilde{A}^r_{k,m}$ are dense. We will do it for $A^r_{n,k,m}$ and leave the case of singularities $\tilde{A}^r_{k,m}$ to the reader. \par

Let $U$ be an open subset of ${\cal X}^r_\omega(M)$ in the $C^r$ topology. We want to show that $U\cap A^r_{n,k,m}\not=\emptyset$.

The set of $C^r$ vector fields that satisfy (P1) is open and dense in ${\cal X}^r_\omega(M)$. \par

Let $X\in U$ be a vector field that satisfies (P1). By the Regularization Lemma we may assume $X$ to be $C^\infty$. Denote by $\theta_1,\dots,\theta_l$ the set of hyperbolic periodic orbits of $X$ of period $\leq n$. Let $D_i$ be a fundamental domain of $W^u(\theta_i)$. Let $V\subset U$ be a neighborhood of $X$ such that any $Y\in V$ satisfies (P1) and has exactly $l$ hyperbolic periodic orbits of period $\leq n$. We are going to make a sequence of perturbations on $X$ to obtain new vector fields, but these will always belong to $V$. Besides that, each vector field will be obtained from the previous one through the Return Lemma, where both vector fields are  $C^\infty$. So all vector fields will be $C^\infty$.\par

For each $1\leq  i\leq l$, choose a finite subset $Z_i$, $1/k$ dense in $D_i^\circ$, and let $Z=\cup_{i=1}^lZ_i$. \par 

We will need the following lemma whose proof we postpone to the end. 
 
\begin{lem} 
Given $m\in \NN$ and a finite subset $Z$ of $\cup_{i=1}^lD^\circ_i$, there is $Y\in V$ such that: \par
	1) $Y_t(D_i)=X_t(D_i)$, for all $t\leq 0$ and $1\leq  i\leq l$.\par
	2) For every $z \in Z$,  $Y_t(z)\in B_{1/m}(z)$ for some $t>m$.
	\label{lem5}
\end{lem}

Let us show how the lemma implies the theorem. Condition 1) of the lemma implies that the local invariant manifolds $\cup_{t\leq 0} X_t(D_i)$ and $\cup_{t\leq0}Y_t(D_i)$ are the same set. Therefore $X$ and $Y$ have the same periodic orbits, and $Y$ satisfies (P1). The sets $D_i^\circ$ and $Z_i$ remain the same for both $X$ and $Y$, implying that $Y$ satisfies  (A1). By 2) $Y$ satisfies (A2), and therefore (P2). So $Y\in U\cap A^r_{n,k,m}$, which proves the theorem in the case of periodic orbits.\par
It only remains to prove the lemma. The argument is on the cardinality of $Z$.\par

Let us start with the case when $Z$ has one element, $Z=\{z\}$. $z\in D_i^\circ$ for some $i$.

We take a small ball $B_\delta(z)$ around $z$ so that $B_\delta(z)\cap\Sigma^u\subset\ D_i^\circ$ and $\delta < {1/m}$.
Now, let $\Sigma$ be the transversal section at $z$ used  to define the Poincaré map, $\Sigma^u$, $D$, etc, as defined in section 2. The Return Lemma says that for the neighborhood $V$ of $X$, the regular point $z$ and the ball $B_\delta(z)$, there exists $X^1\in V$ such that the following holds: if we consider the section $\Sigma \cap B_\delta(z)$  transversal to $X$, there exist $\tau>0$ and $\Sigma_0$, a neighborhood of $z$ in $\Sigma \cap B_\delta(z)$, such that the set $X_{[0,\tau]}(\Sigma_0)=\{X_t(x)|\;\; x\in \Sigma_0, 0<t<\tau\}$ is contained in $B_\delta(z)$, outside this set $X=X^1$, and $X^1_t(z)\in B_\delta(z)$ for some $t>m$.  Therefore $X^1$ satisfies 2). To see that $X^1$ satisfies 1),  note that $X=X^1$ outside $X_{[0,\tau]}(\Sigma_0)$. Then, $X_t^1(x)=X_t(x)$ for $x\in D_i$ and $t\leq0$, implying that $X^1$ satisfies 1).

Now we assume that the lemma holds for sets of cardinality $s$, and write $Z=\{z_1,\dots,z_{s+1}\}$. Each $z_j$ belongs to some $D_i^\circ$ \par
  
By induction there exists $X^2\in V$ that satisfies the following:\par
1') $X^2_t(D_i)=X_t(D_i)$ for all $t\leq 0$ and $1\leq  i\leq l$.\par
2') For $1\leq  j\leq s$,  $X^2_t(z_j)\in B_{\delta_j}(z_j)$ for some $t>m$.

Here the numbers $\delta_j$ are chosen to be less than $1/m$, and small enough so that the balls $B_{\delta_j}(z_j)$ are pairwise disjoint and $B_{\delta_j}(z_j)\cap\Sigma^u\subset\ D_i^\circ$, when $z_j\in D_i^\circ$. We may also assume that the two vector fields coincide outside these balls. 

Let $k$ be such that $z_{s+1}\in D_k^\circ$. Take a ball $B_{\delta_{s+1}}(z_{s+1})$ such that $\delta_{s+1}<{1/m}$, $B_{\delta_{s+1}}(z_{s+1})\cap\Sigma^u\subset\ D_k^\circ$, and the balls $B_{\delta_j}(z_j)$ are pairwise disjoint for $1\leq  j\leq s+1$.

As in step one of the induction, we apply the Return Lemma to obtain a vector field $Y\in V$ such that: \par
 A) $Y\equiv X^2$ outside $B_{\delta_{s+1}}$.\par
 B) $Y_t(D_k)=X^2_t(D_k)$, for all $t\leq 0$.\par
 C) $Y_t(z_{s+1})\in B_{\delta_{s+1}}(z_{s+1})$ for some $t>m$.\par
Here we have to be careful and take $\delta_{s+1}$ small enough so that $B_{\delta_{s+1}}(z_{s+1})$ does not intersect the first return arcs of the points $z_1,z_2,...,z_s$ under $X^2$, and the union of the sets $X^2_t(D_i)$, $1\leq i\leq l$ and $t\leq 0$. Then, $Y$ satisfies items 1 and 2, for $1\leq j\leq s+1$. So this proves the lemma and the proposition, when $r<\infty$.

It is important to make the following remark. We wanted to keep the same sets $Z_i$ for all vector fields so that they would remain $1/k$ dense in $D_i^\circ$ all the way. \par
When $r=\infty$, recall that the $C^\infty$ topology on  ${\cal X}^\infty_\omega(M)$ is defined as the union of the topologies induced by the inclusion maps $i_r: {\cal X}^\infty_\omega(M) \to {\cal X}^r_\omega(M)$, where ${\cal X}^r_\omega(M)$ is equipped with the $C^r$ topology.\par
So, if we define $A^\infty_{n,k,m}$ and $\tilde{A}^\infty_{k,m}$ as the sets of $X \in {\cal X}^\infty_\omega(M)$ that satisfy (P1) , (P2) and (S1), (S2), respectively, then it is immediate that these sets are open and dense in the $C^\infty$ topology, and $${\cal R^+}=\left(\bigcap_{n,k,m}A^\infty_{n,k,m}\right)\cap \left(\bigcap_{k,m}\tilde{A}^\infty_{k,m}\right)$$ is the set required by the proposition.\QED

 
\section{Proof of the Theorem}
In this section we are going to prove the main result of this paper.\\ 

\noindent{\bf Theorem.} {\it Let $M$ be a compact manifold  without boundary and dimension greater or equal to 3. Let $\omega$ be a volume form on $M$ and consider the set ${\cal X}^r_\omega(M)$ of volume preserving vector fields on $M$, with the $C^r$ topology. Then,  ${\cal X}^r_\omega(M)$ contains a residual subset $\cal R$ such that any $X\in {\cal R}$ has the following property: the closure of an invariant manifold of hyperbolic critical elements is a chain transitive set.}\\ 

We will need the following:

\begin{pro}
	 Let  $X$ be a compact connect metric space. If $X$ contains a dense subset of $\omega$-recurrent points, then $X$ is chain transitive. Similarly, if $X$ has a dense subset of $\alpha$-recurrent points, then $X$ is chain transitive.
	\label{lema6}  
	\end{pro}
{\it Proof.} Firstly, suppose that $X$ contains a dense subset of $\omega$-recurrent points. Given any $p,q\in X$, $\epsilon>0$ and $t>0$, we are going to prove  that there exists an $(\epsilon,t)$-chain from $p$ to $q$. Indeed, let $t'>t$. Since $X$ is connect, there exists a finite sequence $(y_0,\dots, y_n)$ in $X$ such that $y_0= \varphi_{t'}(p)$, $y_n=q$ and $d(y_i,y_{i+1})<\epsilon/4$, for all $0\leq i\leq n-1$. For each $0\leq i\leq n-1,$ let $x_i$ be a $\omega$-recurrent point in $B_{\epsilon/2}(y_i)$. Then, it is easy to see that  $[p,x_1,\dots,x_{n-1},q]$ is a $(\epsilon,t)$-chain from $p$ to $q$.

The second part follows from the first. Consider the flow $\psi_t=\varphi_{-t}$. By hypothesis, $X$ has a dense subset of $\omega$-recurrent points, with respect  to $\psi_t$. Then, from the first part,  $X$ is chain transitive with respect  to $\psi_t$. Given $a,b\in X$, $\epsilon>0$ and $t>0$, taking $T>t$ there exists a $(\epsilon,  t)$-chain, with respect to $\psi_t$, from $b$ to $\varphi_T(a)$, say $[x_0=b,x_1,\dots,x_m=\varphi_T(a)]$. For each $0\leq j\leq m-1$, let $t_j>t$ such that $d(\psi_{t_j}(x_j),x_{j+1})<\epsilon$. So, $[a, \psi_{t_{m-1}}(x_{m-1}),\psi_{t_{m-2}}(x_{m-2}),\dots,\psi_{t_{0}}(x_{0}),b]$ is a $(\epsilon,t)$-chain from $a$ to $b$, with respect to $\varphi_t$.  This concludes the proof. \QED\\

The proof follows immediately from Proposition \ref{rec} and Proposition \ref{lema6}, by taking $X$ to be the closure of an invariant manifold.


\begin{thebibliography}{99}
\bibitem{abraham} {\sc Abraham, R.; Marsden, J.} - Foundations of Mechanics. Benjamin/Cummings,  Reading, MA,  1978.
\bibitem{bessa} {\sc Bessa, M.} The Lyapunov exponents of generic zero divergence three-dimensional vector fields. {\it Ergod. Th. \& Dynam. Sys.} {\bf 27} (2007),  1445–1472.
\bibitem{conley} {\sc Conley, C.} - Isolated Invariant Sets and the Morse Index, CBMS Reg. Conf. Ser. in Math., Vol. 38, AMS, Providence, RI (1978).
\bibitem{herman} {\sc Herman, M.} - Examples of compact hypersurfaces in $R^{2p}, 2p\geq 6$, with no periodic orbits. From: ``Hamiltonian systems with three or more degrees of freedom (S'Agar\'o, 1995)", NATO Adv. Sci. Inst. Ser. C Math. Phys. Sci. 533, Kluwer Acad. Publ. (1999), 126. 
\bibitem{oliveira} {\sc Oliveira, F.} -  Density of recurrent points on invariant manifolds of sympletic and volume-preserving diffeomorphisms. {\it Ergod. Th. \& Dynam. Sys.}  {\bf 22} 
(2002), 925-934. 
\bibitem{palis} {\sc Palis, J.; de Melo, W.} - Geometric Theory of Dynamical Systems. Springer-Verlag, 1982. (Translated from Portuguese by
A. K. Manning)
\bibitem{robinson1} {\sc Robinson, C.} - Generic properties of conservative systems. {\it Amer. J. Math.} {\bf 102} (3) (1970), 562-603.
\bibitem{mane} {\sc Ma\~n\'e , R.} - Ergodic Theory and Differentiable Dynamics. Springer-Verlag, 1987. (Translated from Portuguese by
Silvio Levy)
\bibitem{zuppa} {\sc Zuppa, C.} - Regularisation $C^\infty$ des champs vectoriels qui pr\'eservent l'\'el\'ement de volume. {\it Bull. Braz. Math. Soc.} {\bf 10} (2) (1979), 51-56.
\end{thebibliography}
\end{document}